\documentclass[12pt]{article}
\usepackage{amsmath,amssymb}
\usepackage{verbatim}
\usepackage{dsfont,bm}
\usepackage{color}
\usepackage{graphicx}
\usepackage{amsthm}
\usepackage{enumerate}

\usepackage{xcolor}
\definecolor{darkgreen}{rgb}{0,0.75,0}
\definecolor{darkred}{rgb}{0.75,0,0}
\definecolor{darkmagenta}{rgb}{0.5,0,0.5}
\usepackage[colorlinks,citecolor=darkgreen,linkcolor=darkred,urlcolor=darkmagenta]{hyperref}
\newtheorem{theorem}{Theorem}[section]

\newtheorem{lemma}[theorem]{Lemma}

\newtheorem{prop}[theorem]{Proposition}

\newtheorem{remark}[theorem]{Remark}

\numberwithin{equation}{section}

\def\be{\begin{equation}}
\def\ee{\end{equation}}
\def\bes{\begin{equation*}}
\def\ees{\end{equation*}}

\newcommand{\mr}[1]{{\tt \href{http://www.ams.org/mathscinet-getitem?mr=#1}{MR#1}}}
\newcommand{\arxiv}[1]{{\tt \href{http://arxiv.org/abs/#1}{arXiv:#1}}}
\newcommand{\set}[1]{\left\{ #1 \right\}}

\newcommand{\Sett}[2]{\left\{ #1  : \, #2 \right\}}

\newcommand{\abs}[1]{{\left\vert\kern-0.25ex #1
		\kern-0.25ex\right\vert}}
\newcommand\norm[1]{\left\lVert#1\right\rVert} 
\newcommand{\loc}[0]{\operatorname{loc}}

  \def\sC {{\mathcal C}}
 \def\sE {{\mathcal E}} \def\sF {{\mathcal F}}
 \def\sH {{\mathcal H}}

  \def\sX {{\mathcal X}}

 \def\bE {{\mathbb E}} 
\def\bG {{\mathbb G}}  
  
 \def\bN {{\mathbb N}} 
\def\bP {{\mathbb P}}  \def\bR {{\mathbb R}}
  
  \def\bX {{\mathbb X}}
 \def\bZ {{\mathbb Z}}

\def\med{\medbreak\noindent}

\def\sm{\smallskip\noindent}

\def\ignore#1{}

\def\ol{\overline}           

\def\Gam{\Gamma}

\def\to {\rightarrow}

\def\dint{\int\kern-.6em\int}

\newcommand\restr[2]{{
		\left.\kern-\nulldelimiterspace 
		#1 
		\vphantom{\big|} 
		\right|_{#2} 
	}} 
	
	\def\diam{{\mathop{{\rm diam }}}}
	\def\dist{{\mathop {{\rm dist}}}}

	\def\Cap{\operatorname{Cap}}

	\newcommand{\on}[1]{\operatorname{ #1}}
	
	\def\wh{\widehat}
	
	\def\be{\begin{equation}}
	\def\ee{\end{equation}}
	\def\bes{\begin{equation*}}
	\def\ees{\end{equation*}}
	\def\ba{\begin{align}}
	\def\ea{\end{align}}
	\def\xxea{\end{align}}
\def\bas{\begin{align*}}
\def\eas{\end{align*}}

\def\proof{{\smallskip\noindent {\em Proof. }}}
\def\qed{{\hfill $\square$ \bigskip}}

\definecolor{dgreen}{rgb}{0, 0.6, 0.1}
\definecolor{dblue}{rgb}{0, 0.0, 0.6}
\definecolor{vdblue}{rgb}{0,.08, 0.45}
\definecolor{dred}{rgb}{0.7, 0.0, 0.0}
\definecolor{vdblue}{rgb}{0,.08, 0.45}

\definecolor{purple}{rgb}{0.6, 0.0, 0.6}
\definecolor{mytext}{rgb}{0.1, 0.1, 0.1}

\begin{document}
	
	\font\titlefont=cmbx14 scaled\magstep1
	\title{\titlefont A note on heat kernel estimates, resistance bounds and Poincar\'e inequality}  
	\author{		Mathav Murugan\footnote{Research partially supported by NSERC (Canada).} 	}
	\maketitle
	\vspace{-0.5cm}
\begin{abstract}

\vskip.2cm
Sub-Gaussian heat kernel estimates are typical of fractal graphs.
We show that sub-Gaussian estimates on graphs follow from a Poincar\'e inequality, capacity upper bound, and  a slow volume growth condition. An important feature of this work is that we do not assume elliptic Harnack inequality, cutoff Sobolev inequality, or exit time bounds. 
\end{abstract}

\section{Introduction}

\subsection{Statement of the main result}
Let $G=(V,E)$ be an infinite, locally finite connected graph. We say that $x,y \in V$ are neighbours (denoted by $x \sim y$) if $x$ and $y$ are connected by an edge; that is $\set{x,y} \in E$. Let $d_G: V \times V \to \bR$ denote the combinatorial graph metric on $V$. We assume that $G$ is equipped with a weight (or conductance) that is a symmetric nonnegative function $\mu: V \times V \to [0, \infty)$ such that $\mu(x,y)>0$ if and only if $x \sim y$. We call the pair $(G,\mu)$ a \emph{weighted graph}.

Define $\mu_x= \sum_{y \in V} \mu(x,y)$.
 The weight $\mu$ induces a measure on $V$ that we also denote by $\mu$ (note the abuse of notation) defined as $\mu(A)= \sum_{x \in A} \mu_x$. Let $B_G(x,r)$ denote the closed ball and let $V_G(x,r)$ denotes its measure; that is
 \[
 B_G(x,r):= \set{ y \in V: d_G(x,y) \le r }, \hspace{6mm} V_G(x,r):= \mu(B_G(x,r)).
 \]
 In this work, we consider graphs of \emph{polynomial growth}  with volume growth exponent $d_f$; that is, there exists $C_V>1$ such that 
 \be \label{e:vdf} \tag*{$(\on{V}(d_f))$}
C_V^{-1} r^{d_f} \le  V_G(x,r) \le C_V r^{d_f}, \hspace{4mm} \forall x \in V, \hspace{2mm} \forall r \ge 1.
 \ee
For each $x,y \in V$, we set 
\[
p(x,y):=\frac{\mu(x,y)}{\mu_x}.
\]
We say that a weighted graph $(G,\mu)$ satisfies the \ref{e:p0} condition, if there exists $p_0>0$ such that
\be \label{e:p0} \tag*{$(p_0)$}
p(x,y) \ge p_0 \hspace{3mm} \mbox{whenever $x \sim y$}.
\ee
The \ref{e:p0} condition is a discrete analogue of uniform ellipticity.
We will consider discrete time Markov chain $\set{X_n, n \ge 0, \bP^x, x \in V}$, with transition probabilities $p(x,y)$. The chain $(X_n)_{n \ge 0}$ is symmetric with respect to $\mu$, since
\[
p(x,y) \mu_x= \mu(x,y)=\mu(y,x)= p(y,x) \mu_y.
\]
The associated Markov operator $P$ is given by
\[
P f(x):= \bP^x(f(X_1))= \sum_{y \in V} p(x,y) f(y), \hspace{3mm}  \forall f \in \bR^V. 
\]
By the symmetry of $(X_n)_{n \ge 0}$ with respect to $\mu$, $P$ is self-adjoint on $\ell^2(V,\mu)$. 

For $n \in \bZ_+ := \set{0,1,\ldots}$, let $p_n$ denote the transition function of the random walk $X_n$, i.e.
\[
p_n(x,y)= \bP^x(X_n=y).
\]
Let us denote the heat kernel by
\[
h_n(x,y):= \frac{p_n(x,y)}{\mu_y}. 
\]
In this work, we obtain sufficient conditions that are stable under perturbations for sub-Gaussian heat kernel upper and lower bounds: there exists $C,c>0$ such that
\be \label{e:uhk} \tag*{$(\operatorname{UHK}(d_w))$}
 h_n(x,y) \le \frac{C}{V_G(x,n^{1/d_w})} \exp \left[ - \left( \frac{d_G(x,y)^{d_w}}{Cn} \right)^{\frac{1}{d_w-1}} \right], \forall n \ge 1 , \forall x,y \in V,
\ee
and
\be \label{e:lhk} \tag*{$(\operatorname{LHK}(d_w))$}
(h_n+h_{n+1})(x,y) \ge \frac{c}{V_G(x,n^{1/d_w})}  \exp \left[ - \left( \frac{d_G(x,y)^{d_w}}{cn} \right)^{\frac{1}{d_w-1}} \right], \forall n \ge 1 \vee d_G(x,y).
\ee
We denote the conjunction of \ref{e:uhk} and \ref{e:lhk} as $(\on{HK}(d_w))$. 
The parameter $d_w$ is called the \emph{walk dimension} or \emph{escape time exponent}. 
For any $d_f \in [1,\infty)$ and for any $d_w \in [2,d_f+1]$,
Barlow constructs graphs  that satisfy polynomial growth condition \ref{e:vdf} sub-Gaussian estimates \ref{e:uhk} and \ref{e:lhk} -- see  \cite[Theorem 2]{Bar04} and  \cite[Theorem 3.1]{GT}.
Moreover, these are the complete range of $d_f$ and $d_w$ for which \ref{e:vdf}, sub-Gaussian estimates \ref{e:uhk} and \ref{e:lhk} could possibly hold for graphs.


For $f \in \bR^V$, we define
\be 
\sE_G(f,f) = \frac{1}{2} \sum_{x,y \in V } (f(x)-f(y))^2 \mu(x,y). 
\ee
Let $A,B$ be subsets of $V$. We define the \emph{effective conductance}  (or capacity) between $A$ and $B$ as
\be \label{e:defcap}
\on{Cap}_G(A,B):= \inf \set{\sE_G(f,f): f \in \bR^V, \restr{f}{A} \equiv 1, \restr{f}{B} \equiv 0},
\ee
where $\inf \emptyset = + \infty$. By considering $f = 1_A$ or $1_B$, we see that $\on{Cap}(A,B)< \infty$ if $A \cap B = \emptyset$ and if one of $A,B$ is finite.

We say that a weighted graph $(G,\mu)$ satisfies the \emph{capacity upper bound} \ref{e:capub}, if there exists  $C>1$ such that
\be \label{e:capub} \tag*{$(\on{Cap}(d_w)_{\le})$}
\Cap_G(B_G(x,r), B_G(x,2r)^c) \le C \frac{V_G(x,r)}{r^{d_w}}, \hspace{3mm} \forall x \in V, \hspace{2mm} \forall r >0.
\ee
We remark that the constant `$2$' in $B_G(x,2r)^c$ above is essentially arbitrary and can be replaced by any other constant larger than $1$ using the volume doubling property and a covering argument.

We say that a weighted graph satisfies Poincar\'e inequality \ref{e:pi}, if there exists 
$C_P>1$ such that for every ball $B:=B_G(x_0,r), x_0 \in V, r \ge 0$, 
\be  \label{e:pi} \tag{$\on{PI}(d_w)$}
\sum_{x \in B_G(x_0,r)} (f(x)- f_B)^2\mu_x  \le  C_P r^{d_w} \sum_{x,y \in B(x_0,2r), x \sim y}(f(x)-f(y))^2 \mu(x,y),
\ee
where $f_B= \frac{1}{\mu(B_G(x_0,r))} \sum_{x \in B_G(x_0,r)} f(x)\mu_x$. 

The following theorem is the main result of this work.
\begin{theorem}\label{t:main}
Let $(G,\mu)$ be a weighted graph that satisfies the \ref{e:p0} condition,  \ref{e:vdf}, the Poincar\'e inequality \ref{e:pi}, and the capacity upper bound \ref{e:capub} for some $d_f<1+d_w$. Then $(G,\mu)$ satisfies the two sided  sub-Gaussian heat kernel bounds $(\on{HK}(d_w))$.
\end{theorem}
\begin{remark}
	{\rm
\begin{enumerate}[(a)]
	\item We recall that the converse implication $(\on{HK}(d_w)) \Rightarrow (\on{VD})+(\on{PI}(d_w))+ (\on{Cap}(d_w)_\le)$ is satisfied for graphs satisfying the \ref{e:p0} condition \cite[Theorem 3.1]{GT}. Here $(\on{VD})$ denotes the following volume doubling property: there exists $C_D>1$ such that $V_G(x,2r) \le V_G(x,r)$ for all $x \in V,r>0$.

 \item The `Gaussian case' $d_w=2$ is well-known without any restriction on $d_f$. This follows from the work of Delmotte \cite{Del}, which in turn is a discrete analogue of a celebrated theorem of Grigor'yan and Saloff-Coste  \cite{Gri,Sal}.
 
\item  Theorem \ref{t:main} can be viewed as an evidence towards a more general conjecture of Barlow \cite[Remark 3.17(1)]{Bar13}, and of Grigor'yan, Hu and Lau \cite[Conjecture 4.15]{GHL14}, \cite[p. 1495]{GHL15}. This conjecture suggests that Theorem \ref{t:main} should hold without any restriction on $d_f$ and $d_w$. An analogous conjecture has also been made in \cite{GHH} in the context of jump processes. This work provides the first family of examples with $d_w>2$  that are transient ($d_f>d_w$) where the above conjecture can be verified. We remark that the methods of \cite{BCK} can be easily adapted to prove Theorem \ref{t:main} for the `strongly recurrent case' $d_f<d_w$.

\item In the case $d_f<d_w$ the work of Barlow, Coulhon and Kumagai \cite{BCK} provides a satisfactory characterization of $(\on{HK}(d_w))$.  This work can be viewed as progress towards a question raised in Kumagai's ICM survey \cite[Open Problem III]{Kum}, which asks for simpler characterization of $(\on{HK}(d_w))$ for the case $d_f \ge d_w$. For a family of planar graphs satisfying $d_f=d_w$, a different approach that relies on circle packing of planar graphs can be used to prove the above result \cite[Theorem 6.2]{Mur}.

\end{enumerate}
}\end{remark}
 A characterization of $(\on{HK}(d_w))$ that is stable under perturbations was obtained by \cite{BB} for graphs and later extended by \cite{BBK} to metric measure spaces. The characterization is given using a stronger version of \ref{e:capub} known as the cutoff Sobolev inequality  $(\on{CS}(d_w))$.
  This cutoff Sobolev inequality and its variants are a crucial ingredient for the iteration arguments that go back to the works of De Giorgi, Nash and Moser.
 The cutoff Sobolev inequality and its variants have been very useful to obtain stability results for Harnack inequalities and heat kernel estimates for both diffusions and jump processes \cite{AB, BB,BBK, BM, CKW1,CKW2, GHH,GHL15, Lie, MS17,MS}. However, all known proofs of the cutoff Sobolev inequality involve conditions that are \emph{apriori} difficult to obtain (for example, exit time lower bounds). Motivated by these considerations, there is a need to find a simpler alternative to the cutoff Sobolev inequality that is also stable under perturbations as pointed out in the survey \cite[Open problem 4 in p. 38]{Bar03}. 

\subsection{Outline of the proof}
Instead of studying the random walk on a weighted graph, we study the associated cable process. Roughly speaking, the cable process is a diffusion on the corresponding metric graph obtained by placing a unit interval for each edge. In Section \ref{s:pre}, we collect some preliminaries on Dirichlet forms and cable process associated to a graph.

In Section \ref{s:proof}, we provide the proof of Theorem \ref{t:main} which we briefly sketch below.
Our approach is to prove an exit time lower bound on balls 
$B_G(x,r)$ which is at least of the order of  $r^{d_w}$ in the smaller concentric ball, say $B_G(x,r/2)$. It is well known that proving such an exit time lower bound is the main difficulty in obtaining the heat kernel bound $(\on{HK}(d_w))$.

The basic idea behind the constraint $d_f<1+d_w$ is that one-dimensional objects are not negligible when $d_f<1+d_w$ -- see Proposition \ref{p:poin} for a precise formulation. If the exit time of $B(x,r)$ is much less than $r^{d_w}$ at a point close to $x$, by the maximum principle, we have `a tentacle' that is at least one-dimensional in which the exit time is too small. Poincar\'e inequality \ref{e:pi} is used to estimate the capacity of this tentacle from below (see Proposition \ref{p:poin} and \eqref{e:te3}), whereas the capacity upper bound \ref{e:capub} is used to show a competing upper bound on the capacity of the tentacle (see Lemma \ref{l:exit2}(b) and \eqref{e:te4}). These upper and lower bounds on the capacity of the tentacle will contradict each other if the exit time is too small near the center of the ball. This implies the desired exit time lower bound. Once we have the exit time lower bound, we appeal to known results to obtain the two-sided heat kernel bounds.

\section{Preliminaries} \label{s:pre}
In this section, we recall some preliminaries on Dirichlet forms and cable processes.
\subsection{Dirichlet forms}
We recall some standard notions on Dirichlet forms from \cite{FOT}. Let $(\sX,d,m)$ be a locally compact, separable, metric measure space, where $m$ is a Radon measure with full support.
Let $(\sE,\sF)$ be a strongly local, regular Dirichlet form
 on $L^2(\sX,m)$ -- see \cite[Sec. 1.1]{FOT}. Associated with this 
form  $(\sE,\sF)$, there exists an $m$-symmetric Hunt process  $\bX= (\Omega, \sF_\infty,\sF_t,X_t, \bP_x)$ \cite[Theorem 7.2.1]{FOT}.
We denote the \emph{extended Dirichlet space} by $\sF_e$ \cite[Theorem 1.5.2]{FOT}.
 For $f \in \sC_c(\sX) \cap \sF$,  the \emph{energy measure} is defined as the unique Borel measure $d\Gamma(f,f)$ on $\sX$ that satisfies
$$ \int g d\Gam(f,f) = 2 \sE(f,fg) - \sE(f^2,g), \hspace{3mm}\mbox{ for all $g \in \sF \cap \sC_c(\sX)$.} $$
This notion can be extended to all functions in $\sF_e$ and
we have
$$ \sE(f,f) = \int_{\sX} d\Gam(f,f). $$
This follows from \cite[Lemma 3.2.3]{FOT} with a caveat that our definition of $\Gamma(f,f)$ is different from \cite{FOT} by a factor $1/2$.

Let $(\sX,d)$ be a metric space equipped  with a  strongly local Dirichlet form $(\sE,\sF)$ on $L^2(\sX,m)$.
We call  $(\sX,d,m,\sE,\sF)$ a {\em  metric  measure space with Dirichlet form}, or {\em MMD space}.

For an open subset of $D$ of $\sX$, we define the following function spaces 
associated with $(\sE,\sF)$ on $L^2(\sX,m)$.
\begin{align}
\sF_{\loc}(D) &= \Sett{u \in L^2_{\loc}(D,m)}{ \forall \mbox{ relatively compact open } 
	\Omega \subset D, \exists u^\# \in \sF, u = u^\# \big|_{\Omega}  \mbox{ $m$-a.e.}}, \nonumber \\
\sF(D)   &=  \Sett{u \in \sF_{\loc}(D) }{ \int_D \abs{u}^2 \,dm + \int_D d\Gamma(u,u) < \infty }, \label{e:loc1} \\
\label{e:loc2}\sF_{\mbox{\tiny{c}}}(D)  
&=  \Sett{u \in \sF(D)}{\mbox{the essential support of } u \mbox{ is compact in } D }, \\
\sF_0(D)  &= \mbox{ the closure of } \sF_{\mbox{\tiny{c}}}(D) \mbox{ in } \sF \mbox{ in the norm }
\left(\sE(\cdot,\cdot)+ \norm{\cdot}_2^2\right)^{1/2}. \label{e:loc3}
\end{align}

We define \emph{capacities} for a MMD space $(\sX,d,m,\sE,\sF)$  as follows.
By $A \Subset D$, we mean that the closure of $A$ is a compact subset of $D$.
For $A \Subset D$ we set
\be \label{e:capdef}
\Cap_D(A) = \inf\{ \sE(f,f): f \in \sF_0(D) \mbox{ and  $f \ge 1$ in a neighbourhood of $A$} \}.
\ee
The following \emph{domain monotonicity} of capacity is clear from the definition: if $A_1 \subset A_2 
\Subset D_1 \subset D_2$ then
\be \label{e:capmon}
\Cap_{D_2}(A_1) \le \Cap_{D_1}(A_2). 
\ee
Given an open set $\Omega \subset \sX$, a linear operator $G^\Omega: L^2(\Omega) \to \sF$ is called a \emph{Green operator} if, for any $\phi \in \mathcal{C}_c(\Omega) \cap \sF$, and for any $f \in L^2(\Omega)$, 
\be \label{e:gr}
\sE(G^\Omega f, \phi)= \int_{\Omega} f \phi \, d m.
\ee
For an open set $\Omega \subset \sX$, we define $E^\Omega: \Omega \to \bR$ as
\be \label{e:Ex}
E^\Omega:= G^\Omega 1_\Omega. 
\ee
The function $E^\Omega$ has the following probabilistic meaning: $E^\Omega(x)$ is the mean exit time from $\Omega$ of the Hunt process $X_t$ associated to $(\sE,\sF)$ on $L^2(\sX,m)$ started at $x$; that is,
\bes
E^\Omega(x)= \bE_x \tau_\Omega, 
\ees
where
$\tau_\Omega:= \inf \set{t>0: X_t \in \Omega^c}$ is the exit time from $\Omega$. 
\subsection{Cable processes}

In this work, we embed the graph in a connected metric space by replacing each edge by an isometric copy of the unit interval $[0,1]$ and gluing them at endpoints in an obvious manner. This connected metric space is known as the cable system corresponding to the graph.
Random walks on graphs can be studied using a diffusion on the cable system because many relevant properties like Harnack inequalities, heat kernel bounds, exit time estimates, functional inequalities can be transferred between a graph and its cable system \cite{BB}. We refer to \cite{Fol} for an introduction to cable systems and the associated diffusions.

Let $(G,\mu)$ be a weighted graph. We define an associated MMD space $(\sX,d,m,\sE,\sF)$ called the  \emph{cable system} corresponding to $(G,\mu)$. The metric space $(\sX,d)$ is defined as follows. 
We view the edges $E$ as a subset of the two-element subsets of $V$, \emph{i.e.} $E \subset \set{ J \subset V : \abs{J}=2}$.
We define an arbitrary orientation by providing each edge $e \in E$ with a source $\wh{s}:E \to V$ and a target $\wh{t}:E \to V$ such that $e= \set{\wh{s}(e),\wh{t}(e)}$. We say two vertices  $u,v \in V$ are \emph{neighbours} if $\set{u,v} \in E$. We say two distinct edges $e,e' \in E$ are \emph{incident} if $e \cap e' \neq \emptyset$.
We define $\sX=\sX(G,\mu)$ corresponding to the graph $G$ as the topological space obtained
by replacing each edge $e \in E$ by a copy of the unit interval $[0,1]$, glued together in the obvious way, with the endpoints corresponding to the vertices.
More formally, we define $\sX$ as the quotient space $\left(E \times [0,1]\right)/\sim$, where $\sim$ is the smallest equivalence relation such that $\wh{t}(e)=\wh{s}(e')$ implies $(e,1) \sim (e',0)$,  $\wh{s}(e)=\wh{s}(e')$ implies $(e,0) \sim (e',0)$, and  $\wh{t}(e)=\wh{t}(e')$ implies $(e,1) \sim (e',1)$.
Here $E \times [0,1]$ is equipped with the product topology with $E$ being a discrete topological space.
It is easy to check that the topological space above does not depend on the choice of the edge orientations given by $s,t:E \to V$. 
This defines the cable system $\sX$ as a topological space equipped with the canonical quotient map
$q: E \times [0,1] \to \sX$. 
 By \cite[Corollary 3.1.24 and Exercise 3.2.14]{BBI}, there is a unique maximal metric $d:\sX \times \sX \to [0,\infty)$ such that $d(q(e,s),q(e,t)) \le \abs{s-t}$ for all $e \in E_\bG, s,t \in [0,1]$.
This metric space $(\sX,d)$ is also  called a \emph{metric graph} or \emph{one-dimensional polyhedral complex} -- \cite[Section 3.2.2]{BBI}. 

Next, we define the measure $m$ on $\sX$ corresponding to $(G,\mu)$. The  measure $m$ on $(\sX,d)$ is defined as the unique Borel measure that satisfies
\[
m( q(e \times [s,t])) = \mu(u,v) \abs{s-t}, \mbox{ for all $e=\set{u,v} \in E$ and $0 \le s \le t \le 1$.}
\]
We say that $x \in \sX$ is a \emph{vertex} if $x=q(e,0)$ or $x=q(e,1)$ for some $e \in E$. We note that the set of all vertices in $\sX$ have  zero measure.

For an edge $e \in E$, we denote by $\gamma_e: [0,1] \to \sX$ the map $\gamma_e(t)=q(e,t)$. We say that a function $f:\sX \to \bR$ is absolutely continuous if $f$ is continuous and the $f \circ \gamma_e:[0,1] \to \bR$ is absolutely continuous for all $x \in \sX$. If $f$ is absolutely continuous, the function $x \mapsto \abs{\nabla f (x)}$ is well-defined for $m$-almost every $x \in \sX$, by the equation 
\[
\abs{\nabla f}(\gamma_e(s)) = \abs{(f \circ \gamma_e)'(s)}, \hspace{3mm} \mbox{for all $e \in E, s\in [0,1]$.}
\]
Note that although  the sign of $(f \circ \gamma_e)'$ depends on the choice of orientation, the absolute value $\abs{(f \circ \gamma_e)'}$ is independent of the choice of orientation.
Let $\sF_0$ denote the space of compactly supported Lipschitz functions.  By Rademacher theorem, if $f\in \sF_0$, then $f$ is absolutely continuous, $\abs{\nabla f}$ is uniformly bounded and compactly supported, and therefore $\int_{\sX} \abs{\nabla f}^2\,dm < \infty$. We define the  Dirichlet form $(\sE,\sF)$ as
\[
\sE(f,f) = \int_{\sX} \abs{\nabla f}^2(x) \, m(dx),
\]
where $\sF$ is the completion of $\sF_0$ with respect to the norm
\[
\norm{f}_{\sE_1} = \left(\sE(f,f)+ \int_{\sX} f^2 \, dm\right)^{1/2}.
\]
It is easy to verify that the form $(\sE,\sF)$ on $L^2(\sX,m)$ is a closed, Markovian, bilinear form and therefore defines a Dirichlet form \cite{Fol}. In this context, the energy measure is given by
\[
\int_A d\Gamma(f,f) = \int_A \abs{\nabla f}^2 \, dm.
\]

In the remainder of the work, we assume that \ref{e:p0} is \emph{always} satisfied by a weighted graph. For the metric measure space $(\sX,d,m)$, we denote
\[
B(x,r) = \set{y \in \sX: d(x,y)<r}, \hspace{3mm} V(x,r)=m(B(x,r)).
\]

We recall a few properties of the weighted graph $(G,\mu)$ that are inherited by the associated cable system $(\sX,d,m,\sE,\sF)$ \cite[Section 3]{BB}. For instance, the weighted graph $(G,\mu)$ satisfies the polynomial volume  growth condition \ref{e:vdf} if and only if
the cable system $(\sX,d,\mu)$ satisfies the following property: there exists $C>1$ such that
\be \label{e:cvdf}
C^{-1} (r \vee r^{d_f})\le  m(B(x,r)) \le C (r \vee r^{d_f}), \hspace{3mm} \mbox{ for all $x \in \sX, r>0$.}
\ee 
Similarly, the weighted graph $(G,\mu)$ satisfies the Poincar\'e inequality \ref{e:pi}, if and only if the cable system $(\sX,d,\mu,\sE,\sF)$ satisfies the following Poincar\'e inequality \cite[Proposition 3.5]{BB}:
 there exists constants $C,A \ge 1$ such that 
 for all $x\in \sX$, $r \ge 1$ and $f \in \sF$
 \be \label{e:cpi}
 \int_{B(x,r)} (f - \ol f)^2 \,dm \le C (r^2 \vee r^{d_w}) \, \int_{B(x,AR)}d\Gamma(f,f),
 \ee
 where $\ol f= m(B(x,r))^{-1} \int_{B(x,r)} f\, dm$.
 
 Similarly, the weighted graph $(G,\mu)$ satisfies the capacity upper bound \ref{e:capub}, if and only if the cable system $(\sX,d,\mu,\sE,\sF)$ satisfies the following analogous bound \cite[Lemma 2.6]{BB}:
 there exist $C>0$ such that for all $r>0$, $x \in \sX$ 
		\be \label{e:ccapub}
		\Cap_{B(x,2r)}(B(x,r)) \le C_1 \frac{m(B(x,r))}{r^2 \vee r^{d_w}}.
		\ee

\section{Proof of the main result} \label{s:proof}
Our first lemma provides an uniform upper bound on exit times of a ball \eqref{e:exub}, a logarithmic Caccioppoli inequality for the exit time \eqref{e:logcac}, and  an `averaged version' of lower bound on the exit time \eqref{e:exalb}.
\begin{lemma} \label{l:exit2}
Let $(G,\mu)$ be a weighted graph. Let $(\sX,d,m,\sE,\sF)$ denote the corresponding cable system and let $E^\Omega$ denote the exit time of the corresponding cable process as defined in \eqref{e:Ex}.
\begin{enumerate}[(a)]
\item  If  $(G,\mu)$ satisfies \ref{e:vdf} and \ref{e:pi}, then we have the following exit time upper bound: there exists $C>0$ such that
\be
\norm{E^{B(x,r)}}_\infty \le C (r^2 \vee r^{d_w}), \label{e:exub}
\ee
for all $x \in \sX, r>0$.
\item   If  $(G,\mu)$ satisfies \ref{e:capub}, then we have the following bounds: there exists $C>0$ such that for all $x \in \sX, r>0$, we have, writing $u=E^{B(x,r)}$,
 \begin{align}
\int_{B(x,r/2)} \abs{\nabla (\log u)}^2 \, dm &\le C\frac{m(B(x,r))}{r^2 \vee r^{d_w}}, \label{e:logcac}\\
\int_{B(x,r/2)} u^{-1} \, dm &\le C \frac{m(B(x,r))}{r^2 \vee r^{d_w}}. \label{e:exalb}
\end{align}

\end{enumerate}
\end{lemma}
\proof
(a)
The exit time upper bound \eqref{e:exub} follows from \cite[Lemma 4.14]{CKW1} (Although \cite{CKW1} concerns jump processes, the same proof applies to diffusions as well). 

(b)
The proof of \eqref{e:logcac}, \eqref{e:exalb} follows from an argument similar to \cite[Lemma 7.1]{GHL15}. We use the notation $B:=B(x,r)$ and $B/2=B(x,r/2)$ below. Let $\phi \in \sC_c(B) \cap \sF$ be a cut-off function for $B/2 \subset B$ such that $\sE(\phi,\phi) \le 2 \Cap(B/2,B)$. 
\begin{align*}
\int_B \phi^2 \, d\Gamma(\log u,\log u) &= - \int \phi^2 \, d\Gamma(u, u^{-1}) = -\sE(u, \phi^2 u^{-1}) + 2 \int \phi u^{-1} \, d\Gamma(u,\phi) \\\
&\le - \int \phi^2 u^{-1}\,dm + \frac{1}{2} \int \phi^2 u^{-2} \, d\Gamma(u,u) + 2 \int d\Gamma(\phi,\phi)\\
& \le - \int \phi^2 u^{-1}\,dm + \frac{1}{2} \int \phi^2 u^{-2} \, d\Gamma(u,u) + 2 \sE(\phi,\phi).
\end{align*}
Therefore, we have
\be \label{e:ext1}
\frac{1}{2}\int_{B/2} \abs{\nabla (\log u)}^2 \, dm +  \int_{B/2} u^{-1}\,dm  \le \frac{1}{2}\int_B \phi^2 \, d\Gamma(\log u,\log u) +  \int \phi^2 u^{-1}\,dm \le 2 \sE(\phi,\phi).
\ee
Combining \eqref{e:ext1} with the capacity upper bound \ref{e:capub}, we obtain \eqref{e:logcac} and \eqref{e:exalb}.
\qed

\noindent
Recall that the \emph{Hausdorff $s$-content} of a set $E$ in a metric space $(\sX,d)$ is the number
\[
\sH^\infty_s(E) = \inf \sum_{i} r_i^s,
\]
where the infimum is taken over all countable covers of the set $E$ by balls $B_i$ of radius $r_i$. The Hausdorff content enjoys the following monotonicity property:
\be \label{e:moncontent}
\sH^\infty_s (E) \le \sH_s^\infty(F), \hspace{1mm}\mbox{ whenever $E \subseteq F$.}
\ee 
 If $E$ is a connected set in a length space, then the Hausdorff $1$-content is comparable to its diameter as
\be \label{e:diam1content}
\frac{1}{2} \diam(E) \le \sH^\infty_1(E) \le \diam(E).
\ee
The upper bound on $\sH^\infty_1(E)$ is easily obtained by covering $E$ using a single ball while the lower bound is  contained in  \cite[proof of Lemma 2.6.1]{BBI}.

The following lower bound on capacity is an useful consequence of the Poincar\'e inequality.
\begin{prop} \label{p:poin}
Let $(\sX,d,m,\sE,\sF)$ be the cable system corresponding to a weighted graph $(G,\mu)$ that satisfies \ref{e:p0} condition, \ref{e:vdf} and \ref{e:pi}. 
Let $E,F \subset B(x,r)$ be disjoint sets in $\sX$ such that 
$\min(\sH_s^\infty(E),\sH_s^\infty(F)) \ge \lambda r^{s}$, where $d_f \ge s > d_f - d_w$. Then there exists a  constant $C>1$ (that depends only on $s,d_f,d_w, \lambda$ and the constants associated with polynomial volume growth and the Poincar\'e inequality) such that 
\be \label{e:po}
\int_{B(x,9r)} \abs{\nabla u}^2 \, dm \ge C^{-1} r^{d_f-d_w},
\ee
for all $u \in \sF\left(B(x,9r) \right)$ with  $\restr{u}{E} \equiv 1$, $\restr{u}{F} \equiv 0$, where $\sF\left(B(x,9r) \right)$ is as defined in \eqref{e:loc1}.
\end{prop}
\proof
The method of proof goes back to Heinonen and Koskela \cite[Theorem 5.9]{HK} along with some recent ideas in \cite[Theorem 4.5]{Mur}.

By \cite[Lemma 3.1]{Mur}, it suffices to consider the `equilibrium potential' $u(y)=\bP_y( T_E < T_F)$, where $T_E,T_F$ denote the hitting times of sets $E,F$ respectively for the reflected process corresponding to the Dirichlet form $\left(\sE, \sF\left(B(x,9r) \right) \right)$ on $L^2(B(x,9r),m)$.

By applying Poincar\'e inequality (for intervals in $\bR$) on path joining $E$ and $F$ along with $d(E,F) \le 2r$, we obtain that
\be \label{e:po1}
\int_{B(x,9r)} \abs{\nabla u}^2 \, dm \ge \frac{1}{d(E,F)} \ge  \frac{1}{2r}.
\ee
By \eqref{e:po1} along with the fact that $d_w\le  d_f + 1$, we can assume that $r \ge 10$ and $d(E,F) \ge 2$. Henceforth, we shall assume that $r \ge 10$ and $d(E,F) \ge 2$.

 We need the following gradient estimate: there exists $C_2>1$ such that 
 \be \label{e:po2}
 \abs{\nabla u}(y) \le C_2
 \ee
 for almost every $y \in B(x,9r)$. Since $\sX_V$ has measure zero and $\abs{\nabla u }(y) = 0$  for almost every $y \in E\cup F$, it suffices to consider  $y \in B(x,9r) \setminus (\sX_V \cup E \cup F)$.
 
 Every  $y \in B(x,9r) \setminus (\sX_V \cup E \cup F)$  belongs to an unique $\sX_e$ for some edge $e$.
 We consider two cases depending on whether or not $\sX_e \cap (E \cup F)$ is empty.
 If  $\sX_e \cap (E \cup F) = \emptyset$ and both endpoint of $\sX_e$ belongs to $B(x,9r)$, since the value of $u$ at endpoints of $\sX_e$ differ by at most $1$ and $u$ is linear in the edge $\sX_e$, we have $\restr{\abs{\nabla u}}{ \sX_e} \le 1$, which immediately implies \eqref{e:po2}. If  $\sX_e \cap (E \cup F) = \emptyset$ and if one of the end points is not in $B(x,9r)$, then $u$ is constant on $\sX_e$, which implies $\restr{\abs{\nabla u}}{ \sX_e} \equiv 0$.

If $y \in \sX_e \cap (B(x,9r) \setminus (\sX_V \cup E \cup F))$ is such that  $\sX_e \cap (E \cup F) \neq \emptyset$,  then using $\dist_1(E,F) \ge 2$ we have that $\sX_e$ intersects exactly one of the sets $E$ or $F$.
By symmetry, it suffices to consider the case $\sX_e \cap E \neq \emptyset$. Let $I_y$ denote the maximal closed interval on $\sX_e \setminus E$ that contains $y$.
If both the endpoints of $I_y$ belongs to $E$, then it is clear that $u \equiv 1$ on $\sX_e$ and therefore \eqref{e:po2} is satisfied.
Otherwise, consider the vertex $v \in \sX_V \cap I_y$. Consider the cable process starting at the vertex $v$, exiting the star shaped set $I_y \cup \left(  \cup_{e: v \notin \sX_e} \sX_e\right)$.
By the harmonic measure of this star shaped set from \cite[Theorem 2.1]{Fol} and using $u(y)=\bP_y(T_E < T_F)$, we obtain the gradient estimate \eqref{e:po2} in this case as well.

For balls $y \in \sX, s >0$ such that $B(y,s) \subset B(x,9r)$, we define
\be \label{e:po3}
u_{y,s} = \frac{1}{V(y,s)}\int_{B(y,s)} u \,dm.
\ee
By the gradient estimate \eqref{e:po2} and fundamental theorem of calculus, there exists $C_3 >1$ such that for all $y \in \sX, 0<s \le C_3^{-1}$,
we have
\be \label{e:po4}
\abs{u(y)-u_{y,r}} \le  \frac{1}{V(y,r)}\int_{B(y,r)} \abs{u(y) - u(z)} \, m(dz) \le \frac{1}{10}.
\ee

The proof splits into two cases, depending on whether or not there are points $y \in E$ and $z\in F$ so that neither 
\[
\abs{u(y) - u_{y,r}} \hspace{3mm} \mbox{nor} \hspace{3mm} \abs{u(z) -u_{z,2r}}
\] 
exceeds $\frac{1}{5}$. If such points $y \in E, z \in F$ can be found, then
\[
1 \le \abs{u(y)-u(z)} \le \frac{1}{5} + \abs{u_{y,r} - u_{z,2r}} + \frac{1}{5}.
\]
Therefore,  we have
\begin{align*}
\frac{3}{5}\le \abs{u_{y,r} - u_{z,2r}}   &\le \frac{C}{r^{d_f}} \int_{B(z,2r)} \abs{u - u_{z,2r}}\, dm \\
&\le  C r^{-d_f/2} \left( \int_{B(z,2r)} \abs{u - u_{z,2r}}^2\, dm \right)^{1/2}\\
& \le C r^{-d_f/2}  \left( r^{d_w} \int_{B(z,4r)} \abs{\nabla u}^2\, dm \right)^{1/2}  \\
& \le C r^{(d_w-d_f)/2}  \left( \int_{B(x,9r)} \abs{\nabla u}^2\,dm \right)^{1/2}, 
\end{align*}
which  implies \eqref{e:po}. In the above display, we used $B(y,r) \subset B(z,2r)$, and volume growth in the first line, 
Cauchy-Schwarz inequality, and polynomial volume growth in the second line, Poincar\'e inequality \ref{e:cpi} in the third line, and $E \cup F \subset B(x,r)$ along with the triangle inequality in the final line.

The second alternative, by symmetry, is
\be \label{e:po5}
\abs{u(y) - u_{y,r} } \ge \frac{1}{5} \hspace{3mm} \mbox{ for all $y \in E$.}
\ee
Let $i \in \bN$ be the unique integer such that
\be \label{e:po6}
(2 C_3)^{-1}  < 2^{-i} r \le  C_3^{-1}, 
\ee
where $C_3$ is as defined in \eqref{e:po4}.
Using \eqref{e:po4}, and \eqref{e:po5},
we have
\be \label{e:po7}
\abs{u_{y,2^{-i}r} - u_{y,r}} \ge \frac{1}{10} \hspace{3mm} \mbox{ for all $y \in E$}.
\ee
Using Cauchy-Schwarz inequality, Poincar\'e inequality \ref{e:cpi}, \eqref{e:po6} and polynomial volume growth, we obtain the following estimate: for all $y \in E$,
\begin{align*}
1 &\le C  \sum_{j=0}^{i -1} \abs{u_{y,2^{-j}r} - u_{x,2^{-j-1}r}}  \le C  \sum_{j=0}^{i-1} \frac{1}{V(y,2^{-j}r)}\int_{B(y,2^{-j}r)} \abs{u - u_{y,2^{-j}r}} \, dm\\
& \le C  \sum_{j=0}^{i-1} \left( \frac{1}{V(y,2^{-j}r)}\int_{B(y,2^{-j}r)} \abs{u - u_{x,2^{-j}r}}^2 \, dm \right)^{1/2}\\
&\le  C  \sum_{j=0}^{i-1} \left( \frac{(2^{-j}r)^{d_w}}{V(y,2^{-j}r)}\int_{B(y,2^{-j+1}r)} \abs{ \nabla u}^2 \, dm \right)^{1/2} \\
&\le C  \sum_{j=0}^{i-1}  \left((2^{-j}r)^{d_w-d_f}\int_{B(y,2^{-j+1}r)} \abs{ \nabla u}^2 \, dm \right)^{1/2}.
\end{align*}
Therefore, if
\[
\int_{B(y,2^{-j+1}r)} \abs{ \nabla u}^2 \, dm \le \epsilon 2^{-js} r^{d_f-d_w},
\]
for some $\epsilon >0$ and for every $y \in E, 0 \le j \le i -1$, we have that
\[
1 \le C \epsilon^{1/2} \sum_{j=0}^{i-1} 2^{-j(d_w-d_f-s)} \le C \epsilon^{1/2}.
\]
Therefore, for each $y \in E$, there exists an integer $j_y$ with $0 \le j_y \le i-1$, such that
\be \label{e:po8}
\int_{B(y,2^{-j_y+1}r)} \abs{ \nabla u}^2 \, dm \ge \epsilon_0 2^{-j_ys} r^{d_f-d_w},
\ee
for some small enough $\epsilon_0$ depending only on the constants associated with the Poincar\'e inequality and polynomial volume growth. By the $5B$-covering lemma (see \cite[Theorem 1.2]{Hei} or \cite[p. 60]{HKST}), and the separability of $\sX$, there exists a  countable  family of pairwise disjoint balls $B_k=B(y_k,2^{-j_{y_k}+1}r)$, such that
\be \label{e:po9}
E \subset \bigcup_{k} B(y_k,2^{-j_{y_k}+1}5r),
\ee
and, by \eqref{e:po8}, such that
\be \label{e:po10}
\left(\operatorname{radius}\left(B(y_k,2^{-j_{y_k}+1}5r\right) \right)^s\le 2^{-j_{y_k}s+4s}r^s   \le C r^{s+d_w -d_f} \int_{B(y_k,2^{-j_{y_k}+1} r)}  \abs{ \nabla u}^2 \, dm.
\ee
Hence the assumption on $\sH^\infty_s(E)$, \eqref{e:po9}, \eqref{e:po10} and the fact that $B_k$'s are disjoint
\begin{align*}
\lambda r^s &\le \sH_s^\infty(E) \le \sum_k  \left(2^{-j_{y_k}+1} 5 r)\right)^s\\
& \le  C r^{s+d_w-d_f} \sum_{k} \int_{B_k}  \abs{ \nabla u}^2 \, dm \le  C r^{s+d_w-d_f} \int_{B(x,9r)}  \abs{ \nabla u}^2 \, dm.
\end{align*}\qed

The following elementary estimate on $1$-content  will be useful.
\begin{lemma} \label{l:1c}
	Let $(\sX,d,m,\sE,\sF)$ denote the cable system corresponding to weighted graph $(G,\mu)$ that satisfies \ref{e:vdf}.
Then there exists $c_1>0$ that depends only on the constant in \ref{e:vdf} such that
\[
\sH_1^\infty(E) \ge c_1  \left( \left( m(E) \right)^{1/d_f} \wedge m(E) \right)
\]
for any $E \subset \sX$.
\end{lemma}
\proof
Let $B_i=B(x_i,r_i), i \in I$ be countable cover of $E$. Define $I_1= \set{i \in I: r_i \le 1}$,  $I_2= I \setminus I_1$ and 
\[
E_j := \bigcup_{ i \in I_j} B(x_i, r_i), \qquad \mbox{for $j=1,2$.}
\]
We consider two cases depending on whether or not  $m(E_1)$ is greater than $m(E)/2$. 

If $m(E_1) \ge m(E)/2$, then
\[
\sum_{i \in I_1} r_i \ge C_V^{-1} \sum_{i \in I_1} m(B_i) \ge C_V^{-1} m(E_1) \ge C_V^{-1}  m(E)/2.
\]

If $m(E_1)  < m(E)/2$, then $m(E_2) \ge m(E)/2$ and therefore
\[
\sum_{i \in I_2} r_i \ge \left(\sum_{i \in I_2} r_i^{d_f}\right)^{1/d_f} \ge \left(\sum_{i \in I_2} C_V^{-1} m(B_i)\right)^{1/d_f}\ge \left(C_V^{-1}  m(E_2)\right)^{1/d_f}  \ge \left(C_V^{-1}  m(E)/2\right)^{1/d_f}.
\]
In the above display, we use $d_f \ge 1$, and $\operatorname{V}(d_f)$.
Combining the two cases yields the desired result.
\qed

{\sm {\em Proof of Theorem \ref{t:main}.}} 
Let $(\sX,d,m,\sE,\sF)$ be the cable system corresponding to $(G,\mu)$ and let $(X_t)_{t \ge 0}$ denote the associated cable process.

	By \cite[Lemma 2.6]{BB}, \cite[Theorem 1.2]{GHL15} and \eqref{e:exub}, it suffices to verify the following pointwise lower bound on the exit times for $(X_t)_{t \ge 0}$: there exists $c>0$ such that
	for all $x \in \sX, r \ge 1$, and writing $B:=B(x,r)$,
	\be \label{e:te0}
 	\inf_{y \in B(x,r/36)} E^B(y) \ge  c r^{d_w}.
\ee
	Let $u:= E^B=G^B1_B$ denote the exit time function. By \eqref{e:exalb}, Markov's inequality, and $\operatorname{V}(d_f)$, there exists $C_1,C_2>0$ such that $E= \set{y \in B(x,r/18): u(y) \ge C_1^{-1} r^{d_w}}$ satisfies
	\bes 
	m(E) \ge C_2^{-1} r^{d_f}.
	\ees
	Combining this estimate with Lemma \ref{l:1c}, there exists $C_3>1$ which depends only on the constants associated with $\operatorname{V}(d_f)$ such that
		\be \label{e:te1}
		\sH^\infty_1(E) \ge C_3^{-1} r.
		\ee
	Let $K>0$ such that $\inf_{y \in B(x,r/36)} u(y) < e^{-K-1} C_1^{-1} r^{d_w}$. Then by maximum principle, the set 
	$F_K:= \set{y \in  B(x,r/18):  u(y) \le e^{-K-1} C_1^{-1} r^{d_w}}$ contains a path that joins $B(x,r/36)$ and $B(x,(r-\epsilon)/18)^c$ for any $\epsilon >0$.
Therefore by monotonicity of Hausdorff content \eqref{e:moncontent}, and the lower bound of $1$-content  \eqref{e:diam1content}, we have	
	\be \label{e:te2}
	\sH_1^\infty(F_K) \ge \frac{1}{72}r,.
	\ee
	Let $v=\left( \frac{1}{K}\left(1+\log(u) - \log(C_1^{-1} r^{d_w})\right)_+ \right)\wedge 1$. By the logarithmic Caccioppoli  inequality \eqref{e:logcac}, and Markovian property of Dirichlet forms $v \in \sF(B(x,r/2))$, where $\sF(B(x,r/2))$ denotes the domain of Dirichlet form corresponding to the cable process reflected upon hitting $\partial B(x,r/2)$. Using the definitions of $E$ and $F_K$, we can easily check that $\restr{v}{E} \equiv 1$, $\restr{v}{F_K} \equiv 0$.	
Therefore, by Proposition \ref{p:poin} along with \eqref{e:te1}, \eqref{e:te2}, there exists $C_5>0$ such that
\be \label{e:te3}
\int_{B(x,r/2)} \abs{\nabla v}^2 \, dm \ge C_5^{-1} r^{d_f-d_w}.
\ee
	Furthermore, by \eqref{e:logcac} and the contraction property of Dirichlet forms, there exists $C_4>0$ such that
	\be \label{e:te4}
	\int_{B(x,r/2)} \abs{\nabla v}^2 \, dm \le\frac{1}{K^2}\int_{B(x,r/2)} \abs{\nabla (\log u)}^2 \, dm \le \frac{C_4}{K^2} r^{d_f-d_w}.
	\ee
	Combining \eqref{e:te3} and \eqref{e:te4}, we obtain an upper bound $K \le K_0$, where $K_0$ only depends on the constants associated with 	$\operatorname{V}(d_f)$, $\operatorname{PI}(d_w)$, $\operatorname{Cap(d_w)_{\le}}$, $d_f$ and $d_w$. Since our assumption that $F_K$ is non-empty implies an upper bound $K \le K_0$, we obtain \eqref{e:te0} since
	\[
	\inf_{y \in B(x,r/36)} G^B1_B(y) \ge e^{-K_0-2} C_1^{-1} r^d_w, \hspace{3mm} \mbox{ for all $x \in \sX, r \ge 1$.}
	\] 
\qed 
\begin{remark} \label{r:mv} {\rm 
		The proof of Theorem \ref{t:main} implicitly contains an alternate approach to a mean value inequality for superharmonic functions that avoids the usual iteration methods that go back to the works of De Giorgi \cite{DeG} or Moser \cite{Mos}. In the case $d_w>2$, these iteration methods require the cutoff Sobolev inequality \cite{BB,GHL15}. The method of proof used to obtain the exit time lower bound applies verbatim to any non-negative superharmonic function to yield the following mean value inequality under the assumptions of Theorem \ref{t:main}: there exists $0<\theta_1 < \theta_2<1$ and $C>1$ such that for any ball $B(x,r)$ and for any non-negative superharmonic function $u$ in $B(x,r)$, we have
		\[
		 \left(  \frac{1}{m(B(x,\theta_2 r))} \int_{B(x,\theta_2 r)} u^{-1}\, dm\right)^{-1} \le C \inf_{B(x,\theta_1 r)} u.
		\] 
}\end{remark}

\med {\bf Acknowledgement.}
I am grateful to Martin Barlow for providing helpful comments on an earlier draft of this paper. Remark \ref{r:mv} arose due to a question from Moritz Kassmann.  I thank Naotaka Kajino and Takashi Kumagai for their interest in this work.

\noindent Department of Mathematics, University of British Columbia,
Vancouver, BC V6T 1Z2, Canada. \\
mathav@math.ubc.ca

\end{document}